\documentclass[final]{siamltex}
\usepackage{graphicx}
\graphicspath{{Figs/}}

\newcommand\bfA{\mbox{$\mathbf{A}$}}
\newcommand\bfB{\mbox{$\mathbf{B}$}}
\newcommand\bfE{\mbox{$\mathbf{E}$}}
\newcommand\bfH{\mbox{$\mathbf{H}$}}
\newcommand\bfJ{\mbox{$\mathbf{J}$}}
\newcommand\bfU{\mbox{$\mathbf{U}$}}
\newcommand{\half}{{\textstyle{1 \over 2}}}
\newcommand{\eps}{\varepsilon}

\title {Implicit Integration of the Time-Dependent
        Ginzburg--Landau Equations of Superconductivity}
\author{D.~O.~Gunter\and H.~G.~Kaper\and G.~K.~Leaf\thanks
        {Mathematics and Computer Science Division,
         Argonne National Laboratory, Argonne, IL~60439
         (authorname@mcs.anl.gov).
         This work was supported by the
         Mathematical, Information, and
         Computational Sciences Division
         subprogram of Advanced Scientific
         Computing Research,
         U.S. Department of Energy,
         under Contract W-31-109-Eng-38.
        }
       }

\begin{document}
\maketitle

\begin{abstract}
This article is concerned with the integration of the
time-dependent Ginzburg--Landau (TDGL) equations of
superconductivity.
Four algorithms, ranging from fully explicit
to fully implicit, are presented and evaluated
for stability, accuracy, and compute time.
The benchmark problem for the evaluation
is the equilibration of a vortex configuration
in a superconductor that is embedded in a thin insulator
and subject to an applied magnetic field.
\end{abstract}



\pagestyle{myheadings}
\thispagestyle{plain}
\markboth{D.~O.~GUNTER, H.~G.~KAPER, AND G.~K.~LEAF}
         {INTEGRATION OF THE GINZBURG--LANDAU EQUATIONS}

\section{Introduction}
\label{sec-intro}
At the macroscopic level,
the state of a superconductor
can be described in terms of
a complex-valued order parameter
and a real vector potential.
These variables, which determine
the superconducting and electromagnetic
properties of the system at equilibrium,
are found as solutions of the
Ginzburg--Landau (GL) equations
of superconductivity.
They correspond to critical points
of the GL energy functional~\cite{GL,Tinkham-96},
so in principle they can be determined
by minimizing a functional.
In practice, one introduces a time-like
variable and computes equilibrium
states by integrating the time-dependent
Ginzburg--Landau (TDGL) equations.
The TDGL equations, first formulated
by Schmid~\cite{Schmid-66}
and subsequently derived from microscopic
principles by
Gor'kov and \'{E}liashberg~\cite{Gorkov-Eliashberg-68},
are nontrivial generalizations
of the (time-independent) GL equations,
because the time rate of change must be
introduced in such a manner that gauge
invariance is preserved at all times.

We are interested, in particular,
in vortex solutions of the GL equations.
These are singular solutions,
where the phase of the order parameter
changes by $2\pi$ along any closed contour
surrounding a vortex point.
Vortices are of critical importance
in technological applications of
superconductivity.

Computing vortex solutions of the GL equations
by integrating the TDGL equations to equilibrium
has the advantage that the solutions thus found
are stable.
At the same time, one obtains information
about the transient behavior of the system.
Integrating the TDGL equations to equilibrium
is, however, a time-consuming process
requiring considerable computing resources.
In simulations of vortex dynamics
in superconductors, which were
performed on an IBM SP with tens of processors in parallel,
using a simple one-step Euler integration procedure,
we routinely experienced equilibration times
on the order of one hundred
hours~\cite{Braun-96,Crabtree-96,Crabtree-99}.
Incremental changes would gradually drive
the system to lower energy levels.
These very long equilibration times arise,
of course, because we are dealing with
large physical systems undergoing a phase transition.
The energy landscape for such systems is a
broad, gently undulating plain with many
shallow local minima.
It is therefore important to develop
efficient integration techniques that
remain stable and accurate as the
time step increases.

In this article we present four integration
techniques ranging from fully explicit to
fully implicit for problems on
rectangular domains in two dimensions.
These two-dimensional domains should be viewed
as cross sections of three-dimensional
systems that are infinite and homogeneous
in the third direction (orthogonal to the
plane of the cross section), which is the
direction of the field.
The algorithms are scalable in a
multiprocessing environment and
generalize to three dimensions.
We evaluate the performance of each algorithm
on the same benchmark problem, namely,
the equilibration of a vortex configuration
in a system consisting of a superconducting
core embedded in a blanket of insulating material (air)
and undergoing a transition from the Meissner state
to the vortex state under the influence of
an externally applied magnetic field.
We determine the maximum allowable time step
for stability, the number of time steps
needed to reach the equilibrium configuration,
and the CPU cost per time step.

Different algorithms correspond to
different dynamics through state space,
so the eventual equilibrium vortex configuration
may differ from one algorithm to another.
Hence, once we have the equilibrium
configurations, we need some measure
to assess their accuracy.
For this purpose we use three parameters:
the number of vortices,
the mean intervortex distance (bond length),
and the mean bond angle taken over
nearest-neighbor pairs of bonds.
When each of these parameters differs less
than a specified tolerance, we say that
the corresponding vortex configurations
are the same.

Our investigations show that
one can increase the time step
by almost two orders of magnitude,
without losing stability, by going from
the fully explicit to the fully implicit algorithm.
The fully implicit algorithm has a higher cost
per time step, but the wall clock time needed
to compute the equilibrium solution (the most
important measure for practical purposes)
is still significantly less.
All algorithms yield the same
equilibrium vortex configuration.

In Section~\ref{sec-GLM}, we present
the Ginzburg--Landau model of superconductivity,
first in its formulation as a system of
partial differential equations,
then as a system of ordinary differential
equations after the spatial variations have
been approximated by finite differences.
In Section~\ref{sec-alg}, we give
four algorithms to integrate the system
of ordinary equations:
a fully explicit,
a semi-implit,
an implicit, and
a fully implicit algorithm.
In Section~\ref{sec-eval},
we present and evaluate the results
of the investigation.
The conclusions are summarized
in Section~\ref{sec-concl}.

\section{Ginzburg--Landau Model}
\label{sec-GLM}
The time-dependent Ginzburg--Landau (TDGL)
equations of superconductivity~\cite{Tinkham-96,
Schmid-66, Gorkov-Eliashberg-68}
are two coupled partial differential equations
for the complex-valued \textit{order parameter}
$\psi = |\psi| {\rm e}^{i\phi}$
and the real vector-valued \textit{vector potential} $\bfA$,
\begin{eqnarray}
  \label{p-dim}
  \frac{\hbar^2}{2m_sD}
  \left( \frac{\partial}{\partial t} + \frac{ie_s}{\hbar} \Phi \right) \psi
  &=&
  - \frac{1}{2m_s}
  \left( \frac{\hbar}{i} \nabla - \frac{e_s}{c} \bfA \right)^2 \psi
  + a \psi - b |\psi|^2 \psi , \\
  \label{A-dim}
  \nu
  \left( \frac{1}{c} \frac{\partial\bfA}{\partial t} + \nabla \Phi \right)
  &=&
  - \frac{c}{4\pi} \nabla\times \nabla\times \bfA + \bfJ_s .
\end{eqnarray}
Here, $\bfJ_s$ is the \textit{supercurrent density},
which is a nonlinear function of $\psi$ and $\bfA$,
\begin{equation}
  \bfJ_s
  =
  \frac{e_s\hbar}{2im_s} (\psi^* \nabla \psi - \psi \nabla \psi^*)
  - \frac{e_s^2}{m_sc} |\psi|^2 \bfA
  =
  \frac{e_s}{m_s} |\psi|^2
  \left( \hbar \nabla \phi - \frac{e_s}{c} \bfA \right) .
  \label{J-dim}
\end{equation}
The real scalar-valued \textit{electric potential}~$\Phi$
is a diagnostic variable.
The constants in the equations are
$\hbar$,~Planck's constant divided by $2\pi$;
$a$ and $b$, two positive constants;
$c$,~the speed of light;
$m_s$ and $e_s$, the effective mass and charge, respectively,
of the superconducting charge carriers (Cooper pairs);
$\nu$,~the electrical conductivity;
and $D$,~the diffusion coefficient.
As usual,
$i$~is the imaginary unit, and
$^*$~denotes complex conjugation.

The quantity $|\psi|^2$ represents
the local density of Cooper pairs.
The local time rate of change $\partial_t\bfA$
of $\bfA$ determines the \textit{electric field},
$\bfE = (1/c) \partial_t \bfA + \nabla \Phi$,
its spatial variation the (induced) \textit{magnetic field},
$\bfB = \nabla\times \bfA$.

The TDGL equations describe the gradient flow
for the Ginzburg--Landau energy,
which is the sum of
the kinetic energy,
the condensation energy,
and the field energy,
\begin{equation}
  E
  =
  \int
  \left[
  \frac{1}{2m_s}
  \left|
  \left( \frac{\hbar}{i} \nabla - \frac{e_s}{c} \bfA \right)
  \psi
  \right|^2
  +
  \left( - a |\psi|^2 + \frac{b}{2} |\psi|^4 \right)
  +
  | \nabla\times \bfA|^2
  \right]
  \, {\rm d}x .
  \label{E-dim}
\end{equation}
A thermodynamic equilibrium configuration corresponds
to a minimum of $E$.

The energy functional~(\ref{E-dim}) assumes
that there are no defects in the superconductor.
Material defects can be naturally present
or artifically induced and can be in the form
of point, planar, or columnar defects (quenched disorder).
A material defect results in a local reduction of
the depth of the well of the condensation energy.
A simple way to include material defects
in the Ginzburg--Landau model is by assuming
that the parameter $a$ depends on position
and has a smaller value wherever a defect is present.

\subsection{Dimensionless Form}
Let $\psi_\infty^2 = a/b$,
and let $\lambda$, $\xi$, and $H_c$
denote the London penetration depth,
the coherence length, and
the thermodynamic critical field,
respectively,
\begin{equation}
  \lambda = \left( \frac{m_sc^2}{4\pi\psi_\infty^2e_s^2} \right)^{1/2} ,
  \quad
  \xi = \left( \frac{\hbar^2}{2m_sa} \right)^{1/2} ,
  \quad
  H_c = (4\pi a \psi_\infty^2)^{1/2} .
\end{equation}
In this study, we render the TDGL equations dimensionless
by measuring lengths in units of $\xi$,
time in units of the relaxation time $\xi^2/D$,
fields in units of $H_c\surd{2}$,
and energy densities in units of $(1/4\pi)H_c^2$.
The nondimensional TDGL equations are
\begin{eqnarray}
  \label{p}
  \left( \frac{\partial}{\partial t} + i \Phi \right) \psi
  &=&
  \left( \nabla - \frac{i}{\kappa} \bfA \right)^2 \psi
  + \tau \psi - |\psi|^2 \psi , \\
  \label{A}
  \sigma
  \left( \frac{\partial\bfA}{\partial t} + \kappa \nabla \Phi \right)
  &=&
  - \nabla\times \nabla\times \bfA + \bfJ_s ,
\end{eqnarray}
where
\begin{equation}
  \bfJ_s
  =
  \frac{1}{2i\kappa} (\psi^* \nabla \psi - \psi \nabla \psi^*)
  - \frac{1}{\kappa^2} |\psi|^2 \bfA
  =
  \frac{1}{\kappa} |\psi|^2 \left(\nabla \phi -\frac{1}{\kappa} \bfA \right) .
  \label{J}
\end{equation}
Here, $\kappa = \lambda / \xi$
is the Ginzburg--Landau parameter
and $\sigma$ is a dimensionless resistivity,
$\sigma = (4\pi D/c^2) \nu$.
The coefficient $\tau$ has been inserted
to account for defects;
$\tau(x) < 1$ if $x$ is in a defective region;
otherwise $\tau(x) = 1$.
The nondimensional TDGL equations are associated
with the dimensionless energy functional
\begin{equation}
  E
  =
  \int
  \left[
  \left|
  \left( \nabla - \frac{i}{\kappa} \bfA \right)
  \psi
  \right|^2
  +
  \left( - \tau |\psi|^2 + \half |\psi|^4 \right)
  +
  | \nabla\times \bfA|^2
  \right]
  \, {\rm d}x .
  \label{E}
\end{equation}

\subsection{Gauge Choice}
The (nondimensional) TDGL equations are invariant
under a gauge transformation,
\begin{equation}
  {\cal G}_\chi :
  (\psi, \bfA, \Phi)
  \mapsto
  (\psi {\rm e}^{i\chi},
  \bfA + \kappa \nabla \chi,
  \Phi - \partial_t \chi) .
\end{equation}
Here, $\chi$ can be any real scalar-valued function
of position and time.
We maintain the zero-electric potential gauge,
$\Phi = 0$, at all times, using the
\textit{link variable} $\bfU$,
\begin{equation}
  \bfU = \exp \left( -\frac{i}{\kappa} \int \bfA \right) .
  \label{U}
\end{equation}
This definition is componentwise:
$U_x = \exp ( -i\kappa^{-1} \int^x A_x (x',y,z) \, {\rm d}x' )$, $\ldots\,$.
The gauged TDGL equations can now be written
in the form
\begin{eqnarray}
  \label{p-gauge}
  \frac{\partial \psi}{\partial t}
  &=&
  \sum_{\mu = x,y,z}
  U_\mu^* \frac{\partial^2}{\partial \mu^2} (U_\mu \psi)
  + \tau \psi - |\psi|^2 \psi , \\
  \label{A-gauge}
  \sigma
  \frac{\partial\bfA}{\partial t}
  &=&
  - \nabla\times \nabla\times \bfA + \bfJ_s ,
\end{eqnarray}
where
\begin{equation}
  J_{s,\mu}
  =
  \frac{1}{\kappa}
  \mbox{ Im }
  \left[
  (U_\mu \psi)^*
  \frac{\partial}{\partial \mu}
  (U_\mu \psi)
  \right] ,
  \quad \mu = x, y, z .
  \label{J-gauge}
\end{equation}

\subsection{Two-Dimensional Problems}
From here on we restrict the discussion to problems
on a two-dimensional rectangular domain
(coordinates $x$ and $y$), assuming
boundedness in the $x$ direction
and periodicity in the $y$ direction.
The domain represents a superconducting core
surrounded by a blanket of insulating
material (air) or a normal metal.
The order parameter vanishes outside the
superconductor, and no superconducting
charge carriers leave the superconductor.
The whole system is driven by a time-independent
externally applied magnetic field $\bfH$
that is parallel to the $z$ axis,
$\bfH = (0, 0, H)$.
The vector potential and the supercurrent
have two nonzero components,
$\bfA = (A_x, A_y, 0)$ and $\bfJ_s = (J_x, J_y, 0)$,
while the magnetic field has only one nonzero component,
$\bfB = (0, 0, B)$, where
$B = \partial_x A_y - \partial_y A_x$.

\subsection{Spatial Discretization}
The physical configuration to be modeled
(superconductor embedded in blanket material)
is periodic in $y$ and bounded in $x$.
In the $x$ direction, we distinguish three
subdomains: an interior subdomain
occupied by the superconducting material
and two subdomains, one on either side,
occupied by the blanket material.
We take the two blanket layers to be
equally thick, but do not assume
that the problem is symmetric around
the midplane.
(Possible sources of asymmetry are
material defects in the system,
surface currents,
and different field strengths
on the two outer surfaces.)

\begin{figure}[htb]
\begin{center}
\includegraphics[height=1.8in]{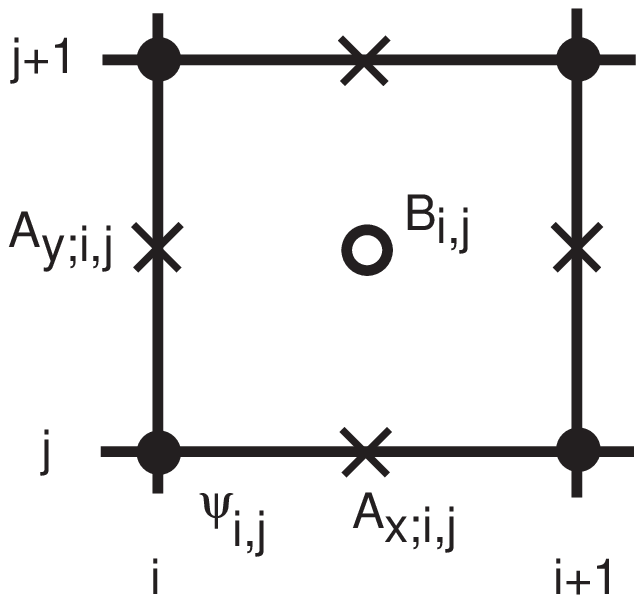}
\end{center}
\caption{Computational cell with evaluation points for
         $\psi$, $A_x$, and $A_y$.}
\label{fig-cell}
\end{figure}

We impose a regular grid
with mesh widths $h_x$ and $h_y$,
\begin{equation}
  \Omega_{i,j} = (x_i, x_{i+1}) \times (y_j, y_{j+1}) , \quad
  x_i = x_0 + i h_x ; \quad
  y_j = y_0 + j h_y ,
  \label{xiyj}
\end{equation}
assuming the following correspondences:
\begin{center}
\begin{tabular}{rll}
  Left outer surface:  & $x = x_0 + \half h_x$,          & $i = 0$,       \\
  Left interface:      & $x = x_{n_{sx}-1} + \half h_x$, & $i = n_{sx}-1$, \\
  Right interface:     & $x = x_{n_{ex}} + \half h_x$,   & $i = n_{ex}$, \\
  Right outer surface: & $x = x_{n_x} + \half h_x$,      & $i = n_x$.
\end{tabular}
\end{center}
\noindent
One period in the $y$ direction is covered
by the points $j = 1, \ldots\,,n_y$.
We use the symbols Sc and Bl to denote the index sets
for the superconducting and blanket region, respectively,
\begin{eqnarray}
  \label{Sc}
  \mbox{Sc} &=& \{ (i,j) :
  (i,j) \in [n_{sx}, n_{ex}] \times [1, n_y] \} , \\
  \label{Bl}
  \mbox{Bl} &=& \{ (i,j) :
  (i,j) \in [1, n_{sx}-1] \cup [n_{ex}+1, n_x] \times [1, n_y] \} .
\end{eqnarray}
The order parameter $\psi$ is evaluated
at the grid vertices,
\begin{equation}
  \psi_{i,j} = \psi (x_i, y_j) , \;
  (i,j) \in \mbox{Sc} ,
  \label{pij}
\end{equation}
the components $A_x$ and $A_y$ of the vector potential
at the midpoints of the respective edges,
\begin{equation}
  A_{x;i,j} = A_x (x_i + \half h_x, y_j) ,\quad
  A_{y;i,j} = A_y (x_i, y_j + \half h_y) , \;
  (i,j) \in \mbox{Sc} \cup \mbox{Bl} ,
  \label{Aij}
\end{equation}
and the induced magnetic field $B$
at the center of a grid cell,
\begin{eqnarray}
  B_{i,j} &=& B (x_i + \half h_x, y_j + \half h_y) \\ \nonumber
  &=& \frac{A_{y;i+1,j} - A_{y;i,j}}{h_x}
  - \frac{A_{x;i,j+1} - A_{x;i,j}}{h_y} , \;
  (i,j) \in \mbox{Sc} \cup \mbox{Bl} ,
  \label{Bij}
\end{eqnarray}
see Fig.~\ref{fig-cell}.
The values of the link variables and the supercurrent
are computed from the expressions
\begin{eqnarray}
  \label{Uij}
  U_{x;i,j}
  = {\rm e}^{-i \kappa^{-1} h_x A_{x;i,j}} &,& \quad
  U_{y;i,j}
  = {\rm e}^{-i \kappa^{-1} h_y A_{y;i,j}} , \\
  \label{Jij}
  J_{x;i,j}
  = \frac{1}{\kappa h_x} {\rm Im}
    \left[ \psi_{i,j}^* U_{x;i,j} \psi_{i+1,j} \right] &,& \quad
  J_{y;i,j}
  = \frac{1}{\kappa h_y} {\rm Im}
    \left[ \psi_{i,j}^* U_{y;i,j} \psi_{i,j+1} \right] .
\end{eqnarray}
The discretized TDGL equations are
\begin{eqnarray}
  \label{p-d}
  \frac{{\rm d}\psi_{i,j}}{{\rm d}t}
  &=&
  \left( L_{xx} (U_{x; \cdot,j}) \psi_{\cdot,j} \right)_i
  + \left( L_{yy} (U_{y; i,\cdot}) \psi_{i,\cdot} \right)_j
  + N \left( \psi_{i,j} \right)
  , \quad
  (i,j) \in \mbox{Sc} ,
  \hspace{-4em}
  \\
  \label{Ax}
  \sigma
  \frac{{\rm d}A_{x;i,j}}{{\rm d} t}
  &=&
  \left( D_{yy} A_{x; i,\cdot}  \right)_j
  - \left( D_{yx} A_{y; \cdot,\cdot}  \right)_{i,j}
  + J_{x;i,j}
  , \quad
  (i,j) \in \mbox{Sc} \cup \mbox{Bl} , \\
  \label{Ay}
  \sigma
  \frac{{\rm d}A_{y;i,j}}{{\rm d} t}
  &=&
  \left( D_{xx} A_{y; \cdot,j}  \right)_i
  - \left( D_{xy} A_{x; \cdot,\cdot}  \right)_{i,j}
  + J_{y;i,j}
  . \quad
  (i,j) \in \mbox{Sc} \cup \mbox{Bl} ,
\end{eqnarray}
where
\begin{eqnarray}
  \hspace*{4em}
  \left( L_{xx} (U_{x; \cdot,j}) \psi_{\cdot,j} \right)_i
  &=& h_x^{-2}
    \left[ U_{x;i,j} \psi_{i+1,j} - 2 \psi_{i,j}
        + U_{x;i-1,j}^* \psi_{i-1,j} \right] , \\
  \left( L_{yy} (U_{y; i,\cdot}) \psi_{i,\cdot} \right)_j
  &=& h_y^{-2}
    \left[ U_{y;i,j} \psi_{i,j+1} - 2 \psi_{i,j}
        + U_{y;i,j-1}^* \psi_{i,j-1} \right] , \\
  N\left( \psi_{i,j} \right)
  &=& \tau_{i,j} \psi_{i,j} - |\psi_{i,j}|^2 \psi_{i,j} , \\
  \left( D_{yy} A_{x; i,\cdot}  \right)_j
  &=& h_y^{-2}
    \left[ A_{x;i,j+1} - 2 A_{x;i,j} + A_{x;i,j-1} \right] , \\
  \left( D_{xx} A_{y; \cdot,j}  \right)_i
  &=& h_x^{-2}
    \left[ A_{y;i+1,j} - 2 A_{y;i,j} + A_{y;i-1,j} \right] , \\
  \left( D_{yx} A_{y; \cdot,\cdot}  \right)_{i,j}
  &=& h_x^{-1}h_y^{-1}
    \left[ \left(A_{y;i+1,j} - A_{y;i,j}\right)
        - \left(A_{y;i+1,j-1} - A_{y;i,j-1} \right) \right] , \\
  \left( D_{xy} A_{x; \cdot,\cdot}  \right)_{i,j}
  &=& h_x^{-1}h_y^{-1}
    \left[ \left(A_{x;i,j+1} - A_{x;i,j}\right)
        - \left(A_{x;i-1,j+1} - A_{x;i-1,j}\right) \right] .
\end{eqnarray}
The interface conditions are
\begin{equation}
  \psi_{n_{sx}-1,j} = U_{x; n_{sx}-1,j} \psi_{n_{sx},j} , \;
  \psi_{n_{ex}+1,j} = U^*_{x;n_{ex},j} \psi_{n_{ex},j} , \;
  j = 1, \ldots\,, n_y .
  \label{p-int}
\end{equation}
At the outer boundary, $B$ is given,
\begin{equation}
  B_{0,j} = H_{L_j} , \;
  B_{n_x,j} = H_{R_j} , \; j = 1, \ldots\,, n_y .
  \label{B-bc}
\end{equation}
The resulting approximation is second-order
accurate~\cite{Gropp-96}.

\section{Time Integration}
\label{sec-alg}
We now address the integration of
Eqs.~(\ref{p-d})--(\ref{Ay}).
The first equation,
which controls the evolution of $\psi$,
involves the second-order linear
finite-difference operators
$L_{xx}$ and $L_{yy}$,
whose coefficients depend on $A_x$ and $A_y$,
and the local nonlinear operator $N$,
which involves neither $A_x$ nor $A_y$.
Each of the other two equations,
which control the evolution of $A_x$ and $A_y$
respectively, involves likewise a second-order
linear finite-difference operator,
but with constant coefficients,
and the nonlinear supercurrent operator,
which involves $\psi$, $A_x$, and $A_y$.
The following algorithms are distinguished
by whether the various operators are treated
explicitly or implicitly.

\subsection{Fully Explicit Integration}
Algorithm~I uses a fully explicit
forward Euler time-marching procedure
for $\psi$, $A_x$, and $A_y$.
Starting from an initial triple
$(\psi^0, A_x^0, A_y^0)$,
we solve for $n = 0, 1, \ldots\,$,
\begin{eqnarray}
  \label{pI-d}
  \frac{\psi^{n+1}_{i,j} - \psi^n_{i,j}}{\Delta t}
  &=& \left( L_{xx} (U^n_{x; \cdot,j}) \psi^n_{\cdot,j} \right)_i
  + \left( L_{yy} (U^n_{y; i,\cdot}) \psi^n_{i,\cdot} \right)_j
  + N \left(\psi^n_{i,j}\right) , \quad
  (i,j) \in \mbox{Sc} ,
  \hspace{-4em}
  \\
  \label{AxI}
  \sigma
  \frac{A^{n+1}_{x;i,j} - A^n_{x;i,j}}{\Delta t}
  &=&
  \left( D_{yy} A^n_{x; i,\cdot}  \right)_j
  - \left( D_{yx} A^n_{y; \cdot,\cdot}  \right)_{i,j}
  + J^n_{x;i,j}
  , \quad
  (i,j) \in \mbox{Sc} \cup \mbox{Bl} , \\
  \label{AyI}
  \sigma
  \frac{A^{n+1}_{y;i,j} - A^n_{y;i,j}}{\Delta t}
  &=&
  \left( D_{xx} A^n_{y; \cdot,j}  \right)_i
  - \left( D_{xy} A^n_{x; \cdot,\cdot}  \right)_{i,j}
  + J^n_{y;i,j}
  . \quad
  (i,j) \in \mbox{Sc} \cup \mbox{Bl} ,
\end{eqnarray}
where $J^n$ is defined in terms of
$\psi^n$, $A^n_x$, and $A^n_y$ in the obvious way.
The initial triple is usually chosen
so the superconductor is in the Meissner state,
with a seed present to trigger the transition
to the vortex state.

Algorithm I has been described in~\cite{Gropp-96}.
It has been implemented in a distributed-memory
multiprocessor environment (IBM SP2);
the transformations necessary to achieve
the parallelism have been described
in~\cite{Galbreath-93}.
The code uses the Message Passing
Interface (MPI) standard~\cite{Dongarra-et-al-98}
as implemented in the MPICH software
library~\cite{Gropp-Lusk-97}
for domain decomposition,
interprocessor communication,
and file I/O.
The code has been used extensively
to study vortex dynamics in superconducting
media~\cite{Braun-96,Crabtree-96,Crabtree-99}.
The underlying algorithm provides
highly accurate solutions
but requires a significant number
of time steps for equilibration.
For stability reasons, the time step
$\Delta t$ cannot exceed 0.0025.

\subsection{Semi-Implicit Integration}
Algorithm~II is generated by an implicit
treatment of the second-order linear
finite-difference operators $D_{yy}$ and $D_{xx}$
in the equations for $A_x$ and $A_y$, respectively,
\begin{eqnarray}
  \label{pII-d}
  \frac{\psi^{n+1}_{i,j} - \psi^n_{i,j}}{\Delta t}
  &=& \left( L_{xx} (U^n_{x; \cdot,j}) \psi^n_{\cdot,j} \right)_i
  + \left( L_{yy} (U^n_{y; i,\cdot}) \psi^n_{i,\cdot} \right)_j
  + N \left(\psi^n_{i,j}\right) , \quad
  (i,j) \in \mbox{Sc} ,
  \hspace{-4em}
  \\
  \label{AxII}
  \sigma
  \frac{A^{n+1}_{x;i,j} - A^n_{x;i,j}}{\Delta t}
  &=&
  \left( D_{yy} A^{n+1}_{x; i,\cdot}  \right)_j
  - \left( D_{yx} A^n_{y; \cdot,\cdot}  \right)_{i,j}
  + J^n_{x;i,j}
  , \quad
  (i,j) \in \mbox{Sc} \cup \mbox{Bl} , \\
  \label{AyII}
  \sigma
  \frac{A^{n+1}_{y;i,j} - A^n_{y;i,j}}{\Delta t}
  &=&
  \left( D_{xx} A^{n+1}_{y; \cdot,j}  \right)_i
  - \left( D_{xy} A^n_{x; \cdot,\cdot}  \right)_{i,j}
  + J^n_{y;i,j}
  . \quad
  (i,j) \in \mbox{Sc} \cup \mbox{Bl} .
\end{eqnarray}
Equations~(\ref{AxII}) and (\ref{AyII})
lead to two linear systems of equations,
\begin{eqnarray}
  \label{PAx}
  \left(
  I - \frac{\Delta t}{\sigma} D_{yy}
  \right) A^{n+1}_{x;i}
  &=&  F_i (\psi^n, A^n_x, A^n_y) ,
  \quad i = 1, \ldots\,, n_x - 1 , \\
  \label{QAy}
  \left(
  I - \frac{\Delta t}{\sigma} D_{xx}
  \right) A^{n+1}_{y;j}
  &=& G_j (\psi^n, A^n_x, A^n_y) ,
  \quad j = 1, \ldots\,, n_y ,
\end{eqnarray}
for the vectors of unknowns
$A_{x;i} = \{ A_{x;i,j} : j = 1, \ldots\,, n_y \}$
and
$A_{y;j} = \{ A_{y;i,j} : i = 1, \ldots\,, n_x-1 \}$.
The matrix $D_{yy}$ has dimension $n_y \times n_y$
and is periodic tridiagonal with elements
$- h_y^{-2}, 2h_y^{-2}, - h_y^{-2}$;
the matrix $D_{xx}$ has dimension $(n_x-1) \times(n_x-1)$
and is tridiagonal with elements
$- h_x^{-2}, 2h_x^{-2}, - h_x^{-2}$,
(except along the edges,
because of the boundary conditions).
Both matrices are independent of $i$ and $j$.
Furthermore, if the boundary conditions
are time independent, they are constant
throughout the time-stepping process.
Hence, the coefficient matrices in
Eqs.~(\ref{PAx}) and (\ref{QAy})
need to be factored only once;
in fact, the factorization can be done
in the preprocessing stage and the factors
can be stored.

In a parallel processing environment,
the coefficient matrices extend
over several processors, so
Eqs.~(\ref{PAx}) and (\ref{QAy})
are broken up in blocks corresponding to the manner in
which the computational mesh is distributed among the
processor set.
We first solve the equations within each
processor (inner iterations)
and then couple the solutions across
processor boundaries (outer iterations).
Hence, we deal with interprocessor coupling
in an iterative fashion.
Two to three inner iterations usually suffice
to reach a desired tolerance for convergence.
After each inner iteration, each processor shares
boundary data with its neighbors through MPI calls.

\subsection{Implicit Integration}
Algorithm~III combines the semi-implicit treatment
of $A_x$ and $A_y$ with an implicit treatment
of the order parameter,
\begin{eqnarray}
  \label{pIII-d}
  \frac{\psi^{n+1}_{i,j} - \psi^n_{i,j}}{\Delta t}
  &=& \left( L_{xx} (U^n_{x; \cdot,j}) \psi^{n+1}_{\cdot,j} \right)_i
  + \left( L_{yy} (U^n_{y; i,\cdot}) \psi^{n+1}_{i,\cdot} \right)_j
  + N \left(\psi^n_{i,j}\right) , \quad
  (i,j) \in \mbox{Sc} ,
  \hspace{-6em}
  \\
  \label{AxIII}
  \sigma
  \frac{A^{n+1}_{x;i,j} - A^n_{x;i,j}}{\Delta t}
  &=&
  \left( D_{yy} A^{n+1}_{x; i,\cdot}  \right)_j
  - \left( D_{yx} A^n_{y; \cdot,\cdot}  \right)_{i,j}
  + J^n_{x;i,j}
  , \quad
  (i,j) \in \mbox{Sc} \cup \mbox{Bl} , \\
  \label{AyIII}
  \sigma
  \frac{A^{n+1}_{y;i,j} - A^n_{y;i,j}}{\Delta t}
  &=&
  \left( D_{xx} A^{n+1}_{y; \cdot,j}  \right)_i
  - \left( D_{xy} A^n_{x; \cdot,\cdot}  \right)_{i,j}
  + J^n_{y;i,j}
  . \quad
  (i,j) \in \mbox{Sc} \cup \mbox{Bl} .
\end{eqnarray}
The second and third equation are solved
as in the semi-implicit algorithm of the
preceding section.
The first equation is solved by a method
similar to the method of Douglas and Gunn~\cite{Douglas-Gunn-64}
for the Laplacian.

We begin by transforming Eq.~(\ref{pIII-d})
into an equation for the correction matrix
$\phi^{n+1} = \psi^{n+1} - \psi^n$.
The equation has the general form
\begin{equation}
  \left( I - \Delta t (L_{xx} + L_{yy}) \right) \phi^{n+1}
  = F(\psi^n, A_x^n, A_y^n) .
  \label{phi}
\end{equation}
If $\Delta t$ is sufficiently small,
we may replace the operator in the left member
by an approximate factorization,
\begin{equation}
  \left( I - \Delta t (L_{xx} + L_{yy}) \right)
  \approx \left( I - \Delta t L_{xx} \right)
  \left( I - \Delta t L_{yy} \right) ,
  \label{fact}
\end{equation}
and consider, instead of Eq.~(\ref{phi}),
\begin{equation}
  \left( I - \Delta t L_{xx} \right)
  \left( I - \Delta t L_{yy} \right)
  \phi^{n+1}
  = F(\psi^n, A_x^n, A_y^n) .
  \label{phi-appr}
\end{equation}
This equation can be solved in two steps,
\begin{eqnarray}
  \label{phi1}
  \left( I - \Delta t L_{xx} \right) \varphi
  &=& F , \\
  \label{phi2}
  \left( I - \Delta t L_{yy} \right)
  \phi^{n+1}
  &=& \varphi .
\end{eqnarray}
The conditions~(\ref{p-int}), which must be satisfied
at the interface between the superconductor and
the blanket material, require some care.
If we impose the conditions at every time step,
then
\begin{eqnarray*}
  \phi^{n+1}_{n_{sx}-1,j}
  &=& U^{n+1}_{x; n_{sx}-1,j}
    \phi^{n+1}_{n_{sx},j}
  + \left[
    U^{n+1}_{x; n_{sx}-1,j} - U^n_{x; n_{sx}-1,j}
    \right] \psi^n_{n_{sx},j} , \\
  \phi^{n+1}_{n_{ex}+1,j}
  &=& \left(U^{n+1}_{x;n_{ex},j}\right)^*
    \phi^{n+1}_{n_{ex},j}
  + \left[
    \left(U^{n+1}_{x; n_{ex},j}\right)^*
    - \left(U^n_{x; n_{sx}-1,j}\right)^*
    \right] \psi^n_{n_{sx},j} ,
\end{eqnarray*}
for $j = 1, \ldots\,, n_y$.
These conditions couple the correction $\phi$
to the update of $A_x$.
To eliminate this coupling, we solve Eq.~(\ref{phi})
subject to the reduced interface conditions
\begin{eqnarray}
  \phi^{n+1}_{n_{sx}-1,j}
  &=& U^{n+1}_{x; n_{sx}-1,j}
    \phi^{n+1}_{n_{sx},j} , \;
  j = 1, \ldots\,, n_y , \\
  \phi^{n+1}_{n_{ex}+1,j}
  &=& \left(U^{n+1}_{x;n_{ex},j}\right)^*
    \phi^{n+1}_{n_{ex},j} ,\;
  j = 1, \ldots\,, n_y .
\end{eqnarray}
When Eq.~(\ref{phi}) is replaced by
Eq.~(\ref{phi-appr}), these conditions are
inherited by the system~(\ref{phi1}).

\subsection{Fully Implicit Integration}
Algorithm~IV uses a fully implicit
integration procedure for the order parameter,
\begin{eqnarray}
  \label{pIV-d}
  \frac{\psi^{n+1}_{i,j} - \psi^n_{i,j}}{\Delta t}
  &=& \left( L_{xx} (U^n_{x; \cdot,j}) \psi^{n+1}_{\cdot,j} \right)_i
  + \left( L_{yy} (U^n_{y; i,\cdot}) \psi^{n+1}_{i,\cdot} \right)_j
  + N \left(\psi^{n+1}_{i,j}\right) , \quad
  (i,j) \in \mbox{Sc} ,
  \hspace{-8em}
  \\
  \label{AxIV}
  \sigma
  \frac{A^{n+1}_{x;i,j} - A^n_{x;i,j}}{\Delta t}
  &=&
  \left( D_{yy} A^{n+1}_{x; i,\cdot}  \right)_j
  - \left( D_{yx} A^n_{y; \cdot,\cdot}  \right)_{i,j}
  + J^n_{x;i,j}
  , \quad
  (i,j) \in \mbox{Sc} \cup \mbox{Bl} , \\
  \label{AyIV}
  \sigma
  \frac{A^{n+1}_{y;i,j} - A^n_{y;i,j}}{\Delta t}
  &=&
  \left( D_{xx} A^{n+1}_{y; \cdot,j}  \right)_i
  - \left( D_{xy} A^n_{x; \cdot,\cdot}  \right)_{i,j}
  + J^n_{y;i,j}
  . \quad
  (i,j) \in \mbox{Sc} \cup \mbox{Bl} .
\end{eqnarray}
The new element here is the term
$N \left(\psi^{n+1}_{i,j}\right)$
in the first equation.

The second and third equations are solved again
as in the semi-implicit algorithm.
The first equation is solved by a slight
modification of the method used in the
implicit algorithm of the preceding section,
The modification is brought about
by the approximation
\begin{equation}
  N \left(\psi^{n+1}\right)
  = \tau \psi^{n+1} - |\psi^{n+1}|^2 \psi^{n+1}
  \approx
  \frac{1}{\Delta t}
  \left( S \left(\psi^n\right) - \psi^n \right) ,
  \label{sg}
\end{equation}
where $S$ is a nonlinear map,
\begin{equation}
  S(\psi)
  = \frac{\tau^{1/2} \psi}
  {\left[
  |\psi|^2 + (\tau - |\psi|^2) \exp(-2\tau \Delta t)
  \right]^{1/2} } .
  \label{S}
\end{equation}
(This approximation is explained in the remark below.)
Equation~(\ref{pIV-d}) is
again of the form~(\ref{phi}),
but with a different right-hand side,
\begin{equation}
  (I - \Delta t (L_{xx} + L_{yy})) \phi^{n+1}
  = G (\psi^n, A_x^n, A_y^n) .
  \label{phi-IV}
\end{equation}
The difference is that,
where $F$ in Eq.~(\ref{phi}) contains a term
$(\Delta t) N \left(\psi^n\right)$,
$G$ in Eq.~(\ref{phi-IV}) contains
the more complicated term
$S \left(\psi^n\right) - \psi^n$.

\paragraph{Remark}
The approximation~(\ref{sg}) is suggested
by semigroup theory. Symbolically,
\begin{equation}
  N(\psi)
  = \lim_{\Delta t \to 0} \frac{S(\Delta t) \psi - \psi}{\Delta t} .
\end{equation}
To find an expression for the ``semigroup'' $S$,
we start from the continuous TDGL equations
(\ref{p})--(\ref{J})
(zero-electric potential gauge, $\Phi = 0$),
using the polar representation
$\psi = |\psi| {\rm e}^{i\phi}$,
\begin{eqnarray}
  \label{TDGL-psi-0}
  \partial_t |\psi|
  &=& \Delta |\psi| - |\psi| | \nabla\phi - \kappa^{-1} \bfA |^2
  + \tau |\psi| - |\psi|^3 , \\
  \label{TDGL-phi-0}
  |\psi| \partial_t \phi
  &=& 2 (\nabla |\psi|) \cdot (\nabla\phi - \kappa^{-1} \bfA)
  + |\psi| \nabla\cdot (\nabla\phi - \kappa^{-1} \bfA) , \\
  \label{TDGL-A-0}
  \sigma \partial_t \bfA
  &=& - \nabla\times \nabla\times \bfA
  + \kappa^{-1} |\psi|^2 (\nabla\phi - \kappa^{-1} \bfA) .
\end{eqnarray}
At this point, we are interested in the effect
of the nonlinear term  $|\psi|^3$ on the dynamics.
To highlight this effect,
we concentrate on the time evolution of
the scalar $u = |\psi|$ and
the vector $v = \nabla\phi - \kappa^{-1} \bfA$.
(In physical terms, $u^2$ is the density of
superconducting charge carriers, while
$u^2v$ is $\kappa$ times the supercurrent density.)
Ignoring their spatial variations, we have
a dynamical system,
\begin{eqnarray}
  \label{TDGL-u-1}
  u' &=& - u | v |^2 + \tau u - u^3 , \\
  \label{TDGL-v-1}
  v' &=& -\eps u^2 v ,
\end{eqnarray}
where ${}'$ denotes differentiation
with respect to $t$, and
$\eps = (\kappa^2\sigma)^{-1}$.
This system yields a pair of
ordinary differential equations
for the scalars $x=u^2$ and $y = |v|^2$,
\begin{eqnarray}
  x' &=& 2 x (\tau - x - y) , \\
  y' &=& - 2 \eps xy .
\end{eqnarray}
If $\kappa$ is large, $\eps$ is small,
and the dynamics are readily analyzed.
To leading order, $y$ is constant;
$y=0$ is the only meaningful choice.
(Recall that $xy^{1/2}$ is $\kappa$ times
the magnitude of the supercurrent density.)
Then the dynamics of $x$ are given by
\begin{equation}
  x' = 2x (\tau - x) .
\end{equation}
We integrate this equation from $t = t_n$ to $t$,
\begin{equation}
  x(t)
  =
  \frac{\tau x(t_n)}{x(t_n) + (\tau - x(t_n))\exp(-2\tau(t - t_n))} .
\end{equation}
In particular,
\begin{equation}
  x(t_{n+1})
  =
  \frac{\tau x(t_n)}{x(t_n) + (\tau - x(t_n))\exp(-2\tau \Delta t)} ,
\end{equation}
where $\Delta t = t_{n+1} - t_n$.
Since $x(t_n) = |\psi^n|^{1/2}$
and $x(t_{n+1}) = |\psi^{n+1}|^{1/2}$,
it follows that
\begin{equation}
  |\psi^{n+1}|
  =
  \frac{\tau^{1/2} |\psi^n|}
       {[|\psi^n|^2 + (\tau - |\psi^n|^2)\exp(-2\tau \Delta t)]^{1/2}} .
\end{equation}
The phase $\phi$ of $\psi$ is constant in time.
If we multiply both sides by ${\rm e}^{i\phi}$, we obtain
the expression~(\ref{S}) for the ``semigroup'' $S$.

\section{Evaluation}
\label{sec-eval}
We now present the results of several
experiments, where the algorithms
described in the preceding section
were applied to a benchmark problem.

\subsection{Benchmark Problem}
The benchmark problem adopted for this
investigation was the equilibration of
a vortex configuration in a superconductor
(Ginzburg-Landau parameter $\kappa = 16$)
embedded in a thin insulator (air),
where the entire system was periodic
in the direction of the free surfaces ($y$).

The superconductor measured $128 \xi$
in the transverse ($x$) direction.
The thickness of the insulating layer on
either side was taken to be $2\xi$,
so the total width of the system was $132 \xi$.
The period in the $y$ direction was taken to be $192 \xi$,
so the entire configuration measured $132 \xi \times 192 \xi$.

The computational grid was uniform,
with a mesh width $h_x = h_y = \half\xi$.
The periodic boundary conditions in the $y$ direction
were handled through ghost points,
so the computational grid had $264 \times 386$ vertices.
The index sets for the superconductor and blanket
(see Eqs.~(\ref{Sc}) and~(\ref{Bl})) were
\begin{eqnarray}
  \mbox{Sc} &=& \{ (i,j) : i = 5, \ldots \,, 260,
                  \, j = 1, \ldots \,, 386 \} , \\
  \mbox{Bl} &=& \{ (i,j) : i = 1, \ldots \,, 4, 261, \ldots \,, 264,\,
                  \, j = 1, \ldots \,, 386 \} .
\end{eqnarray}
The applied field was uniform in $y$ and
equally strong on the left and right side
of the system,
\begin{equation}
  H_L = H_R = H = 0.5 .
\end{equation}
(Units of $H$ are $H_c \surd 2$,
so $H \approx 0.707\ldots H_c$).
As there is no transport current in the system,
the solution of the TDGL equations tends to
an equilibrium state.

\subsection{Benchmark Solution}
\begin{figure}[htb]
\begin{center}
\mbox{\includegraphics[height=3.5in]{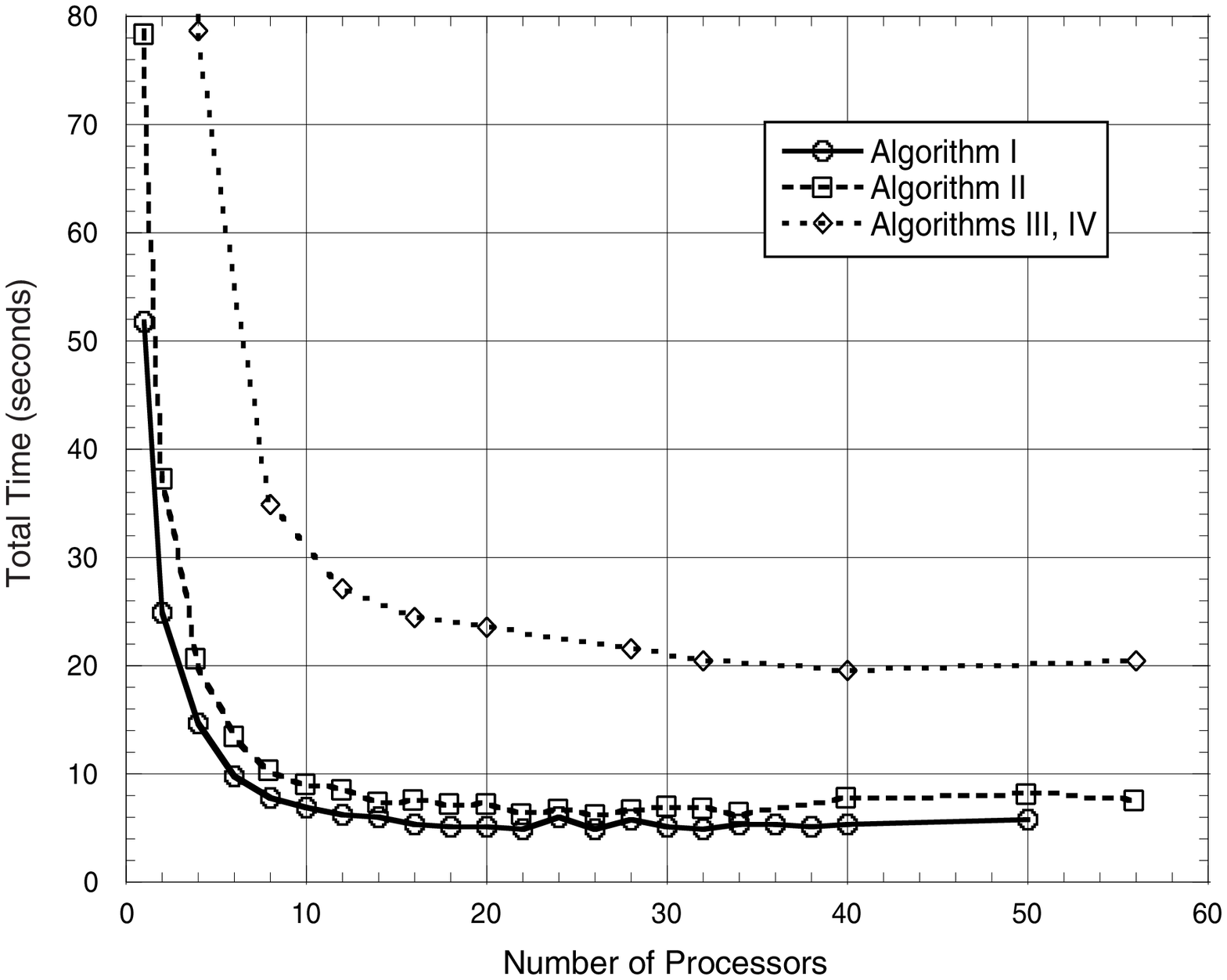}}
\end{center}
\caption{Elapsed time for 50 time steps
         as a function of the number of processors.}
\label{fig-saturation}
\end{figure}

First, preliminary runs were made
to determine, for each algorithm,
the optimal number of processors
in a multiprocessing environment.
Figure~\ref{fig-saturation} shows the elapsed
(wall clock) time for 50 time steps
against the number of processors
on the IBM SP2.
Each algorithm showed a saturation
around 16 processors,
beyond which any improvement became marginal.
All problems were subsequently run on 16 processors.

\begin{figure}[htb]
\begin{center}
\mbox{\includegraphics[height=4.5in]{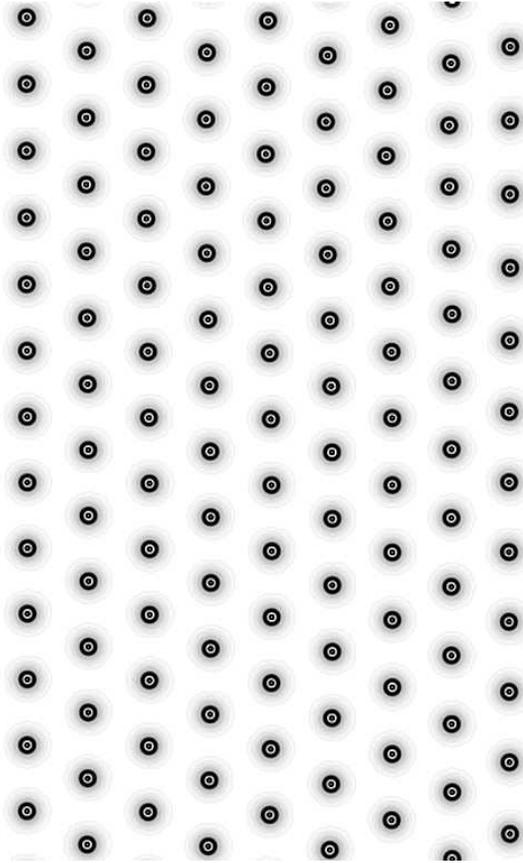}}
\end{center}
\caption{Equilibrium vortex configuration for the benchmark problem.}
\label{fig-config}
\end{figure}

Next, the fully explicit Algorithm~I was used
to establish a benchmark equilibrium configuration.
Equations~(\ref{pI-d})--(\ref{AyI})
were integrated with a time step
$\Delta t = 0.0025$ (units of $\xi^2/D$),
the maximal value for which the algorithm
remained stable.
The evolution of the vortex configuration
was followed by monitoring the number of vortices
and their positions.
Equilibrium was reached after 10,000,000 time steps,
when the number of vortices remained constant
and the vortex positions varied less than
$1.0 \times 10^{-6}$ (units of $\xi$).
The equilibrium vortex configuration had 116 vortices
arranged in a hexagonal pattern; see Fig.~\ref{fig-config}.
The elapsed time for the entire computation
was 50.81 hours.
The elapsed time per time step (0.018 seconds)
is a measure for the computational cost of Algorithm~I.

\subsection{Evaluation of Algorithms II--IV}
Once the benchmark solution was in place,
each of the remaining algorithms (II--IV)
was evaluated for stability, accuracy, and
computational cost.

The stability limit was found by
gradually increasing the time step and
integrating until equilibrium.
Above the stability limit, the algorithm
failed because of arithmetic divergences.
Equilibrium was defined by the same criteria
as for the benchmark solution:
no change in the number of vortices
and a variation in the vortex positions of
less than $1.0 \times 10^{-6}$.
The results are given in Table~\ref{table-summary};
$\Delta t$ is the time step at
the stability limit (units of $\xi^2/D$),
$N$ the number of time steps needed
to reach equilibrium,
$T$ the elapsed (wall clock) time (in hours)
needed to compute the equilibrium configuration,
and $C$ the cost (in seconds per time step,
$C = 3600 T/N$).

Because each algorithm defines its own path through
phase space, one cannot expect to find identical equilibrium
configurations nor equilibrium configurations
that are exactly the same as the benchmark.
The equilibrium vortex configurations
for the four algorithms were indeed different,
albeit slightly.
To measure the differences quantitatively,
we computed the following three parameters:
(1)~the number of vortices in the superconducting region,
(2)~the mean bond length joining neighboring pairs
    of vortices, and
(3)~the mean bond angle subtended by neighboring bonds
    throughout the vortex lattice.
In all cases, the number of vortices was the same (116);
the mean bond length varied less than $1.0 \times 10^{-3} \xi$,
and the mean bond angle varied by less than $1.0 \times 10^{-3}$ radians.
Within these tolerances, the equilibrium vortex
configurations were the same.

\begin{table}[htb]
\begin{center}
\caption{Performance data for Algorithms I--IV.
         \label{table-summary}}

\vspace{2ex}
\begin{tabular}{|| c ||r|r|r|r||}\hline
 &
\multicolumn{1}{|c}{$\Delta t$} &
\multicolumn{1}{|c}{$N$} &
\multicolumn{1}{|c}{$C$} &
\multicolumn{1}{|c||}{$T$} \\
Algorithm    & & & & \\ \hline
I            & 0.0025 & 10,000,000 & 0.018 & 50.81 \\
II           & 0.0500 &    500,000 & 0.103 & 14.32 \\
III          & 0.1000 &    250,000 & 0.232 & 16.11 \\
IV           & 0.1900 &    131,580 & 0.233 &  8.41 \\ \hline
\end{tabular}
\end{center}
\end{table}

\begin{figure}[htb]
\begin{center}
\mbox{\includegraphics[height=1.8in]{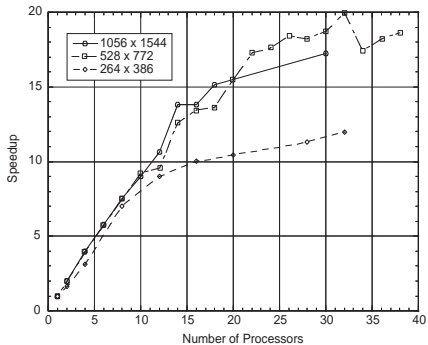}}
\end{center}
\caption{Computational cell with evaluation points for
         $\psi$, $A_x$, and $A_y$.}
\label{fig-speedup}
\end{figure}

Finally, we evaluated the fully implicit Algorithm~IV
from the point of view of parallelism.
From the benchmark problem we derived two more
problems by twice doubling the size of the system
in each direction,
while keeping the mesh width the same ($\half \xi$).
The resulting computational grid had
$528 \times 772$ vertices for the intermediate problem
and $1056 \times 1544$ vertices for the largest problem.
Speedup was defined as the ratio of
the wall clock time (exclusive of I/O)
to reach equilibrium on $p$ processors
divided by the time to reach equilibrium
on a single processor for the smallest
and intermediate problem,
or twice the time to reach equilibrium
on two processors for the largest problem.
(The largest problem did not fit on a single processor.)
The results are given in Fig.~\ref{fig-speedup}.
The curve for the benchmark problem was obtained
as an average over many runs; the data for the
intermediate and largest problem were obtained
from single runs, hence they are less smooth.
The speedup is clearly linear when the number
of processors is small;
it becomes sublinear at about
12 processors for the smallest problem,
14 processors for the intermediate problem,
and 18 processors for the largest problem.

\section{Conclusions}
\label{sec-concl}
The results of the investigation lead
to the following conclusions.

(1)~One can increase the time step $\Delta t$
nearly 80-fold, without losing stability,
by going from the fully explicit Algorithm~I
to the fully implicit Algorithm~IV.

(2)~As one goes to the fully implicit Algorithm~IV,
the complexity of the matrix calculations and,
hence, the cost $C$ of a single time step increase.

(3)~The increase in the cost $C$ per time step
is more than offset by the increase in the size
of the time step $\Delta t$.
In fact, the wall clock time needed to compute
the same equilibrium state with the fully implicit
Algorithm~IV is one-sixth of the wall clock time
for the fully explicit Algorithm~I.

(4)~The (physical) time to reach equilibrium---that is,
$N \Delta t$, the number of time steps needed to reach
equilibrium times the step size---is (approximately)
the same for all algorithms, namely, 25,000
(units of $\xi^2/D$).

(5)~The fully implicit Algorithm~IV displays
linear speedup in a multiprocessing environment.
The speedup curves show sublinear behavior
when the number of processors is large.

\end{document}